\documentclass{article}
\usepackage[utf8]{inputenc}
\usepackage{amsmath,amsfonts,amssymb}
\usepackage{amsthm}
\usepackage[noadjust]{cite}
\newcommand{\e}{\mathsf{E}}
\newcommand{\p}{\mathsf{P}}
\newcommand{\di}{\mathsf{Var}}
\newcommand{\eps}{\varepsilon}

\newtheorem*{cond1}{Condition}
\newtheorem{thm1}{Theorem}
\newtheorem{cor1}{Corollary}
\newtheorem{lm1}{Lemma}

\title{Strong Gaussian approximation for cumulative processes}
\author{Elena Bashtova\footnote{Lomonosov Moscow State University; supported by RFBR grant 20-01-00487; e-mail: elena.bashtova@math.msu.ru}, Alexey Shashkin }
\date{}

\begin{document}

\maketitle

\section{Introduction}

Cumulative processes were initially introduced by Smith in 1955 \cite{Smith1955}. Idea of Smith’s construction of regenerative and cumulative processes was to create an opportunity to apply methods of analysis of sums of independent random variables to a wider class of sequences and processes. Smith proved the stability theorem for regenerative processes and LLN and CLT for cumulative processes. The concept of regeneration turned to be very popular in the queueing theory and was renovated after Smith’s work in a lot of ways. Markov-modulated processes, batch Markov-arrival processes   \cite{Asmussen,Neuts}, regenerative flows in \cite{LG}, cumulative processes as defined in \cite{GW} represent sub-classes of cumulative processes, introduced by Smith. Special definitions arose in order to find explicit formulae for characteristics of models. But as was shown already by Smith himself and later in studies \cite{LG,Thorisson} some classes of important theorems, like  stability conditions or diffusion approximations, may be proved under general assumptions. 
Our conditions {\bf{(A)}} and {\bf{(B)}} below, taken together, are somewhat more restrictive than Smith's definition as we are tackling an almost sure approximation, but rather general and include the most part of examples in literature, e.g. those mentioned above. 

Strong Gaussian approximation of random processes presents an important class of limit theorems (almost sure invariance principles), stemming from seminal results by Skorokhod and Strassen \cite{Strass}. While classical functional limit theorems establish the proximity of properly rescaled  random processes formed by partial sums to  standard ones (e.g. to the Wiener process) in distribution, strong functional limit theorems provide the ways to estimate the almost sure difference of those two types of random processes. For example, given a stationary sequence $ \{X_1,X_2,\ldots\}$ of centered random variables and denoting $S_n=X_1+\ldots+X_n$ for $n\geq 1,$ a strong invariance principle (almost sure invariance principle) would typically state that the sequence can be initially defined on a rich enough probability space so that for some Wiener process $\{W_t,t
\geq 0\}$ and some specified $\sigma^2>0$ (meaning the asymptotic variance of $S_n$) one has
\begin{equation}\label{intr1}
S_n - \sigma W_n = O(f(n)) \mbox{ a.s.,}
\end{equation}
$f=\{f(n),n\geq 1\}$ being some non-random increasing sequence in $n\geq 1$. For \eqref{intr1} to be non-trivial (e.g. allowing to deduce the law of the iterated logarithm for $S_n$) $r$ should satisfy $f(n)=o(\sqrt{n\log\log n})$ which was Strassen's statement. Strong Gaussian approximation allows to deduce almost sure limit theorems, such as laws of the iterated logarithm, by using corresponding statements for approximating Gaussian processes. It has been shown that strong invariance principles  have important applications   to establishment of Darling-Erd\H{o}s type theorems and integral tests for upper and lower functions \cite{StrassBerk,DieEin}.    

One can go further and analyze the rate of convergence in \eqref{intr1}, e.g. obtain an estimate of order $f(n)=n^{1/2-\eps}$ with some nonrandom $\eps>0$.  
  Fundamental results by Koml\'{o}s, Major and Tusn\'{a}dy \cite{KMT,KMT2,Maj} state that the rates $f(n)=\log n$ for i.i.d. sequences with finite exponential moments and $f(n)=n^{1/p}$ if they possess only a moment of order $p>2$ are true and non-improvable. These results were achieved by developing the quantile transformation method. These authors give also a more delicate estimate of the rate of Gaussian approximation in terms of probabilistic maximal inequalities. Namely, in the case of finite exponential moments one has
  \begin{equation}\label{intr2}
\p\left(\sup_{k\leq n}|S_k - \sigma W_k| \geq C\log n+x\right)\leq Ae^{-Bx} 
\end{equation}
for some non-random $A,B,C>0$ and all $n\geq 1,x\geq 0,$ and similar bounds exist in case of lower moment restrictions. Such bounds imply the convergence results of type \eqref{intr1}.

The power and beauty of almost sure (strong) invariance principles prompted extensive research for such results for many classes of random systems. Vector-valued random processes were studied by many authors (\cite{GZ} and references therein), in particular, Zaitsev \cite{Zait} extended \eqref{intr2}  to sums of independent random vectors.  Berkes and Morrow  \cite{BerMor} initiated the exploration of strong approximation for random fields, establishing the rate $f(n)=(n_1\ldots n_d)^{1/2-\eps}$ (here partial sums are indexed by multiindices $n=(n_1,\ldots,n_d)\in\mathbb{N}^d$). For renewal processes, exact rates of convergence in the sense of Koml\'{o}s-Major-Tusn\'{a}dy were established by Cs\"{o}rg\H{o}, 
Horv\'{a}th and Steinebach \cite{CHS} and are in agreement with what one has in the i.i.d. setup (in this case, moment conditions are imposed on the innovations). However for systems of dependent random variables, while there have been many results proved (for martingales \cite{StrassBerk,Eberl,Mermart}, associated random systems  \cite{BSh},  mixing families \cite{PhiStout,Shao,Rio,BerLiuWu} and others), one can rarely attain the optimal convergence rate.  In particular, in most cases methods based on the Skorokhod embedding and the approximation of partial sums by those of independent families do not provide the convergence faster that $O(n^{1/4}),$ which is not optimal in case moments greater than fourth exist. However this border was broken in the works by Berkes, Liu and Wu \cite{BerLiuWu} for processes driven by regular enough Bernoulli shifts, by Merlev${\rm\grave{e}}$de and Rio \cite{MerRio} for functionals of Markov chains and by Cuny, Dedecker, Korepanov, Merlev${\rm\grave{e}}$de \cite{CDKM} for mixing dynamical systems.

  
  This paper is devoted to optimal logarithmic rates of convergence in almost sure invariance principle for multivariate cumulative processes. 
  Our goal is to obtain exponential probabilistic inequalities  of Koml\'{o}s-Major-Tusn\'{a}dy type \eqref{intr2}. To establish the strong Gaussian approximation exact rate under exponential moment restrictions  we follow the approach of Merlev${\rm\grave{e}}$de and Rio \cite{MerRio} to geometrically ergodic Markov chains. Their approach is based on the approximation by an auxiliary process formed by a composition of a Wiener and Poisson ones (Lemma 2.4 \cite{MerRio}). Besides,  the Zaitsev  strong approximation technique \cite{Zait} is applied. Merlev${\rm\grave{e}}$de and Rio study the sequence of partial sums $S_n=\sum_{i=1}^n f(\xi_j),$ $n\geq 1,$ with a bounded Borel functional $f$ and $\{\xi_n,n\geq 1\}$ being a geometrically ergodic Markov chain. However their Lemma 2.4 is a powerful tool of its independent interest.   In particular, we were able to apply it together with Zaitsev construction to analysis of cumulative processes. Note that random sums $S_n$ studied in \cite{MerRio} present a partial case of processes under consideration, as can be seen by taking the Markov chain returns to zero as points of regeneration. 
  Besides, we allow vector-valued processes and dispense with the boundedness condition on the random summands. In the concluding part of the paper we give two corollaries about two important subclasses of cumulative processes. The first one deals with stopped sums. Cs\"{o}rg\H{o}, Deheuvels and Horv\'{a}th \cite{Deheu} obtained the optimal rate in strong invariance principle for such sums  in case when the indexing process is independent from summands, while Horv\'{a}th \cite{Hor} studied general case and got the convergence rate of order $O(n^{1/4}(\log n)^{1/2}(\log\log n)^{1/4})$. We consider stopped sums as in \cite{Hor} and get the logarithmic rate of convergence as a consequence of our main theorem. The second corollary concerns birth and death processes arising in queueing theory \cite{KarlinMc}.
  
 In this paper
 we study a $\mathbb{R}^d-$valued random process $S=\{S(t),t\geq 0\} = \{(S_1(t),\ldots,S_d(t)),t\geq 0\}$ which is assumed to be separable. For a vector $z\in\mathbb{R}^d$, by $|z|$ we denote the maximal norm. The principal condition we will require is: 
 
 \begin{cond1}[{\textbf{A}}]
 There exists an a.s. increasing random sequence $\{T_k, k\geq 0\}$ such that $T_0=0$ and random elements
 $$ \Bigl\{\bigl(T_j-T_{j-1}, S(T_{j-1}+t) - S(T_{j-1}), t\in(0,T_j-T_{j-1}]  \bigr) ,j\geq 1\Bigr\}$$
 are i.i.d.
 \end{cond1}
 
 We need the following notation. 
\begin{itemize} 
\item $\tau_k = T_k - T_{k-1},\;k\geq 1;$
\item $\xi_k =(\xi_{k1},\ldots,\xi_{kd})= S(T_k) - S(T_{k-1}),\; k\geq 1;$
\item  $m(t) = \max\limits\{k\in\mathbb{Z}_+: T_k \leq t\},\; t\geq 0$, i.e. $m=\{m(t),t\geq 0\}$ is the renewal process built by sequence $\{\tau_k,k\geq 1\}$;
\item $\eta_k = \max\limits_{0< t \leq T_k-T_{k-1}}|S(T_{k-1}+t) - S(T_{k-1})|,\;k\geq 1;$ 
\item $\mu=\e \tau_1$, $\varkappa = \mu^{-1}\e \xi_1  ;$
\item $\sigma^2 =   \mu^{-1}(\di(\xi_1) - 2 cov(\xi_1,\tau_1)\varkappa^T + \varkappa \varkappa^T \di(\tau_1)),$
 \end{itemize}
 where $cov(\xi_1,\tau_1) := (cov(\xi_{11},\tau_1),\ldots, cov(\xi_{1d},\tau_1) ).$ We will write $\sigma$ for the usual square root of the matrix $\sigma^2.$

 Next condition restricts the tails of $\eta$ and $\tau$ and is a typical requirement to ensure the logarithmic rate of convergence in the strong invariance principle.
  \begin{cond1}[{\textbf{B}}]
There is $s>0$ such that $\e \exp\{s \tau_1\}<\infty$ and $\e\exp\{s \eta_1\}<\infty$.
 \end{cond1}
 Condition $\mathbf{(B)}$ implies that the exponential moments  of $\xi_1$ are finite in some neighbourhood of the origin.

\section{Main result}

\begin{thm1}\label{TheoremMain} 
Suppose that $\{S(t),t\geq 0\}$ satisfies conditions {\bf(A)} and {\bf (B)}.  Then one can redefine it on a probability space $(\Omega,\mathcal{F},\mathsf{P})$ supporting also a standard $d$-dimensional Wiener process $\{W_t,t\geq 0\}$ such that for some positive constants $a,b,c,$ any $x>0$ and  $t\geq 1$ one has
\begin{equation}\label{1}
\p\left(\sup_{u\leq t}\bigl|S(u) - \varkappa u - \sigma W_u\bigr| \geq c\log t + x \right) \leq a e^{-b x} . 
\end{equation}
\end{thm1}
Standard estimate shows that \eqref{1} implies the strong invariance principle, i.e. the relation
$$S(t) - \varkappa t - \sigma W_t = O(\log t)\mbox{ when }t\to\infty,\mbox{ a.s}. $$

\begin{proof} Without loss of generality we assume  $\di (\tau_1)>0$ as otherwise the statement would be easily derived from the Koml\'{o}s-Major-Tusn\'{a}dy theorem for i.i.d. summands \cite[Theorem 1]{KMT2}. 
We will use triplets of positive constants $(a_i,b_i,c_i),\,i\geq 1,$ when we want to write an auxiliary estimate of the type \eqref{1}.  

Note that we always can confine ourselves to the case $t\geq M,$ with a fixed number $M>1$. For if theorem is proved for that case, then given $t\in[1,M),$ we will be able to write
$$\p\left(\sup_{u\leq t}\bigl|S(u) - \varkappa u - \sigma W_u\bigr| \geq c\log t + x \right) \leq \p\left(\sup_{u\leq M}\bigl|S(u) - \varkappa u - \sigma W_u\bigr| \geq c  + x +c(\log t -\log M) \right)\leq  $$
$$\leq \begin{cases}ae^{-b(x-c\log M)},& x\geq c\log M, \\e^{c\log M-x},& 0< x \leq c\log M.\end{cases} $$
Using this remark we will consider only $t\geq e$ in what follows.

Write
\begin{equation}\label{greeks}
\beta = \bigl(\di (\tau_1)\bigr)^{-1} {cov(\xi_1, \tau_1)},\; v^2 = \di(\xi_1 - \beta \tau_1),\;
\gamma = \frac{\di(\tau_1)}{\mu}, \;\lambda = \frac{\mu^2 }{\di(\tau_1)},\;\alpha = \beta - \varkappa.
\end{equation}

Using the remark above, we will consider only $t\geq e.$
Note that the choice of $\beta$ implies that $$cov(\xi_{ki} - \beta_i \tau_k,\tau_k)=0, \;i=1,\ldots,d. $$


 
Applying Zaitsev  \cite[Theorem 1.3]{Zait} to a $(d+1)$-dimensional sequence $\{(\xi_k - \beta \tau_k +\alpha\mu,\tau_k-\mu), k\geq 1\}$ we can infer (redefining the initial process on a larger probability space if necessary) that there exist two independent Wiener processes $\{B_t,t\geq 0\},\{\widetilde{B}_t,t\geq 0\} $, the former being $d-$dimensional, such that for any pair $(t,x)\in [e,\infty)\times [0,\infty)$

\begin{equation}\label{2}
\p\left(\sup_{k\leq t}\bigl|S(T_k) - k\varkappa \mu - \beta(T_k - k\mu) - vB_k \bigr| \geq c_1\log t + x \right) \leq a_1 e^{-b_1 x} , 
\end{equation}
\begin{equation}\label{3}
\p\left(\sup_{k\leq t}\bigl|T_k - k\mu  - \sqrt{\di(\tau_1)}\widetilde{B}_k \bigr| \geq c_1\log t + x \right) \leq a_1 e^{-b_1 x} , 
\end{equation}
with some $a_1,b_1,c_1>0.$

Arguing as before  (2.50) in Merlev${\rm\grave{e}}$de and Rio \cite{MerRio}, we can construct also a Poisson process $\{N(t),t\geq 0\}$ with parameter $\lambda$, measurable with respect to the $\sigma-$algebra generated by $\widetilde{B}$ and such that the following relations hold:

\begin{equation}\label{4a}
\p\left(\sup_{u\leq t}\bigl|\gamma N_u - u\mu- \sqrt{\di(\tau_1)}\widetilde{B}_u \bigr| \geq c_2\log t + x \right) \leq a_2 e^{-b_2 x} . 
\end{equation}

\begin{equation}\label{5}
\p\left(\sup_{u\leq t}\bigl|\gamma N_u - T_{[u]} \bigr| \geq c_2\log t + x \right) \leq a_2 e^{-b_2 x}  
\end{equation}
with some $a_2,b_2,c_2>0.$

Further, we denote now $y(u) = N^{-1}(u/\gamma)$. Then, process $m$ is the inverse of process $T$ and $y$ is the inverse of   $\gamma N$. Therefore in view of the paper of  Cs\"{o}rg\H{o}, Horv\'{a}th and Steinebach \cite{CHS} and particular Corollary 4.2 there, and of \eqref{3}--\eqref{4a}, one can construct (enlarging the probability space if necessary) a standard Wiener process $\widetilde{W}=\{\widetilde{W}_t,t\geq 0\}$ such that the following statements are simultaneously true:

\begin{equation}
\label{deheuvels}
\p\left(\sup_{u\leq t}\left|m(u) - \frac{u}{\mu} - \frac{1}{\lambda\sqrt{\gamma}}\widetilde{W}_u\right|\geq c_3\log t + x\right)  \leq a_3 e^{-b_3 x},
\end{equation}
 
\begin{equation}\label{6}
\p\left(\sup_{u\leq t}\left| y(u)-  \frac{u}{\mu} - \frac{1}{\lambda\sqrt{\gamma}}\widetilde{W}_u \right| \geq c_3\log t + x \right) \leq a_3 e^{-b_3 x} 
\end{equation}
and
\begin{equation}\label{m-y}
\p\left(\sup_{u\leq t}  |m(u) - {y(u)}| \geq c_3\log t + 
  x \right) \leq a_3e^{  - b_3 x},
 \end{equation}
 with some $a_3,b_3,c_3>0.$ Moreover, process $\widetilde{W}$ is measurable with respect to a $\sigma$-algebra generated by $\widetilde{B}$ plus some random element $V$ independent from all ones considered above.
 
 \begin{lm1}\label{estimate_for_m}
 For any $z\geq 0$ and $t\geq e$ one has
 $$\p\left(\left|m(t) - \frac{t}{\mu} - \frac{1}{\lambda\sqrt{\gamma}}\widetilde{W}_t\right|> z\right)  \leq \max\{a_3,1\} \exp\Bigl\{\max\{b_3,1\}c_3\log t   -\min\{b_3,1\}z \Bigr\}. $$
 \end{lm1}
 
 \begin{proof}
 Consider separately the cases $z\geq c_3\log t$ and $z<c_3\log t,$ and use \eqref{deheuvels} for the first one.
 \end{proof}
 
\begin{lm1}\label{Lemma24} One can construct a standard $d$-dimensional Wiener process $W^*=\{W^*_t,t\geq 0\}$ such that  
\begin{equation}\label{7}\p\left(\sup_{k\leq t}\left|  B_k - \frac{W^*_{N(k)}}{\sqrt{\lambda}}   \right| \geq c_4\log t + x \right) \leq a_4 e^{-b_4 x}  \end{equation}
with some $a_4,b_4,c_4>0$ and all $t\geq e,\,x> 0;$ this process is determined by the process $B$ plus some random element $V^*$ which is independent from all other random elements considered above. In particular, $W^*$  is independent from $\widetilde{B}$ and $N$.
\end{lm1}

\begin{proof} This is a multi-dimensional version of Lemma 2.4 in \cite{MerRio} and is proved along the lines of that lemma, with the following change. One introduces mutually independent Gaussian random vectors
$$Y_{j,k} =2B_{(k+1/2)2^j} -B_{(k+1)2^j}-B_{k2^j},\,j\in\mathbb{Z},k\geq 0 $$
and random functions
$$a_{j,k}=N((k+1/2)2^j) -N(k2^j),\;b_{j,k}=N((k+1)2^j) -N((k+1/2)2^j),$$
$$f_{j,k}=\frac{b_{j,k}I\{(N(k2^j) ,N((k+1/2)2^j)]\}   -a_{j,k} I\{(N((k+1/2)2^j) ,N((k+1)2^j)]\} }{a_{j,k}b_{j,k}(a_{j,k}+b_{j,k})} ,\,j\in\mathbb{Z},k\geq 0$$
where the fraction with $0$ at the denominator is treated as zero. Then, using that $\{f_{j,k}\}_{j,k}$ conditionally on $N$ is an orthonormal family plus some set of identically zero functions, one shows that random vectors
$$W_l=\sum_{j\in\mathbb{Z}}\sum_{k=0}^{\infty}Y_{j,k}\int_0^l f_{j,k}(t)dt ,\,l\geq 0, $$
are jointly Gaussian and $\e W_lW_m^T = \min\{l,m\}I_d $ for $l,m\geq 0,$ where $I_d$ is the identity matrix. The rest of the proof goes like in \cite{MerRio}.
\end{proof}
 
 By construction,  standard Wiener processes $\{\widetilde{W}_t,t\geq 0\}, \{W^*_t,t\geq 0\}$, the latter being $d$-dimensional,  are independent.

\begin{lm1}\label{nonint} In  $\eqref{7}$ one can change the supremum to be over all real $u\in[0,t]$ $($changing $a_4,b_4,c_4$ if necessary$).$
\end{lm1}

\begin{proof}
For a standard scalar Wiener process $Z=\{Z_t,t\geq 0\}$ one has
$$\p\left(\sup_{u\leq t}|Z_u-Z_{[u]}| \geq  \log t +x \right)\leq (t+1)\p\left(\sup_{u\leq 1}|Z_u| \geq \log t + x\right)\leq $$
\begin{equation}\label{Levy}\leq 4(t+1)\exp\left\{-\log^2t/2 - x\log t \right\} \leq    10e^{-x} \end{equation}
as $t\geq e.$
Finally, if $Z$ is independent from the Poisson process $N,$ for $x\geq 1$ we may write
$$\p\left(\sup_{u\leq t}|Z_{N(u)}-Z_{N([u])}| \geq  \log t +x \right)
\leq \sum_{k=0}^{[t]}\p\left( \sup_{u\leq 1}|Z_{N(k+u)}-Z_{N(k)}| \geq  \log t +x \right) \leq $$ $$ \leq (t+1)\p\left(N(1) >\frac{1}{2}(\log t + x)\right)+  (t+1)\p\left( \sup_{u\leq (\log t+x)/2}|Z_{u}| \geq  \log t +x \right) \leq  $$
$$\leq \frac{t+1}{t}e^{-x}\e e^{2N(1)} + (t+1)\p\left( \sup_{u\leq 1}|Z_{u}| \geq  \sqrt{2(\log t +x)} \right)\leq 2e^{\lambda (e^2-1)-x} + 8e^{-x}.  $$
\end{proof}

Let $W^{\circ}=\{W^{\circ}_t,t\geq 0\}$ be a standard $d$-dimensional Wiener process independent from $(W^*,\widetilde{W}).$ Define a Gaussian   process
\begin{equation}\label{4}
W_t = \sigma^{+}\left( \lambda^{-1/2}{v} W^*_{t/\gamma}  - \lambda^{-1}\gamma^{-1/2} \mu \alpha  \widetilde{W}_{t} \right) + (I_d - \sigma^{+} \sigma)W^{\circ}_t, \end{equation}
where $v$ is the square root of the covariance matrix $v^2,$ and $\sigma^+$ stands for the Moore-Penrose pseudo-inverse matrix to $\sigma.$

\begin{lm1}\label{Wiener}
$W$ is a standard $d$-dimensional Wiener process.
\end{lm1}

\begin{proof} It suffices to compute the covariance matrix of $W_1,$ which is equal to
$$\frac{1}{\mu}\sigma^{+}\Bigl(\di(\xi_1-\beta\tau_1) + \di(\tau_1(\beta-\varkappa))\Bigr)\sigma^{+} + I_d - \sigma^{+} \sigma ,$$
and since the components of the random vectors whose covariance matrices are inside the brackets are mutually uncorrelated, using the definition of $\sigma^2$ we infer that $\di W_1=I_d.$
\end{proof}

Denoting $y = y(u) = N^{-1}(u/\gamma) $  we write, for $u\geq 0,$
$$S(u) - \varkappa u -\sigma W_u = S(u)   - \varkappa u - \sigma\sigma^+\frac{v}{\sqrt{\lambda}}W^*_{u/\gamma} + \sigma\sigma^+\alpha \frac{\mu}{\lambda\sqrt{\gamma}}\widetilde{W}_{u} =$$
$$=S(u)   - \varkappa u - \frac{v}{\sqrt{\lambda}}W^*_{u/\gamma}  -\alpha \mu \left(N^{-1}(u/\gamma) - \frac{u}{\lambda \gamma} - \frac{\widetilde{W}_{u }}{\lambda\sqrt{\gamma}}\right) + \alpha\mu\left(N^{-1}(u/\gamma) - \frac{u}{\lambda \gamma} \right) =$$
$$= (S(u) - S(T_{m(u)}) )+   (S(T_{m(u)}) - S(T_{[y]})) +\left(S(T_{[y]}) -\beta T_{[y]}   + \alpha \mu y - v B_y\right)+ $$ $$+\beta  (T_{[y]} - \gamma N(y)) - \alpha \mu \left(y - \frac{u}{\lambda \gamma} - \frac{\widetilde{W}_{u }}{\lambda\sqrt{\gamma}}\right)  + $$ $$+v\left(B_y - \lambda^{-1/2}W^*_{N(y)}\right)+\lambda^{-1/2}v\left(  W^*_{N(y)}-  W^*_{u/\gamma}\right)+\beta(\gamma[u/\gamma]+\gamma-u)  =:\sum_{q=1}^8 \Phi_q(u).  $$
We took into account that $\sigma\sigma^+ v = v $ and  $\sigma\sigma^+ \alpha=\alpha. $
Pick $s_1\in(0,s)$ and  $r>0$ such that \begin{equation}\label{select_r} \e \exp\{s_1|\xi_1 - \varkappa \tau_1|\}<\infty,\; r\log \e \exp\{s_1|\xi_1 - \varkappa \tau_1|\} < \frac{s_1}{6}.\end{equation}

From now till Lemma \ref{LD}, we always consider pairs $(t,x)$ such that 
\begin{equation}\label{pair}x \leq \frac{3 t}{\mu r}.\end{equation}
\begin{lm1}\label{Poisson}
We have $$\p\left(N^{-1}(t/\gamma)\geq 2t/\mu\right) \leq 2e^{  - \lambda (1-2/e)  r x/3},$$
for $(t,x)$ satisfying $\eqref{pair}.$
\end{lm1}

\begin{proof} $N^{-1}(t/\gamma)$ is the sum of $[t/\gamma]+1$ i.i.d. exponential random variables with parameter $\lambda,$ hence
$$\p(N^{-1}(  t/\gamma)\geq 2t/\mu)\leq  e^{-t/\gamma} 2^{[t/\gamma]+1},$$
 and it remains to use the inequality \eqref{pair} and the equality $\mu=\lambda\gamma$.
\end{proof}
   
 \begin{lm1}\label{MerlevedeRio}
There are  $a_5,b_5,c_5>0$ such that for $q=3,\ldots,8$ we have $$\p\left(\sup_{u\leq t}   |\Phi_q(u)|I\{N^{-1}(t/\gamma)\leq 2t/\mu \} \geq c_5\log t + x \right) \leq a_5e^{  - b_5 x},$$
 for all $(t,x)$ satisfying $\eqref{pair}$.
\end{lm1}

\begin{proof}
For $q=3$ one should apply \eqref{2} and \eqref{Levy}.
For $q=4,5,6$ the corresponding statements directly follow from \eqref{5}, \eqref{6}, \eqref{7} respectively with the help of Lemma \ref{nonint}.
For $ q=7$, using that $|u/\gamma - N(y)|\leq 1,$ one has
$$\p\left(\sup_{u\leq t}|W^*_{N(y)}-  W^*_{u/\gamma}|>3d\log t+x\right)\leq
(t/\gamma +1)\p\Bigl(\sup_{u\leq 1}|B_u|>\log t + x/3d\Bigr)\leq $$ $$\leq  4(t/\gamma +1)\exp\{-\log^2 t -(x\log t)/3d\}   \leq 4(\gamma^{-1}+e^{-1})e^{1/2-x/3d}. $$
Finally, for $q=8$ the statement is trivial as $\Phi_8$ is uniformly bounded.
\end{proof}

 \begin{lm1}\label{Phi1}
There are  $a_6,b_6,c_6>0$ such that   $$\p\left(\sup_{u\leq t}   |\Phi_1(u)| \geq c_6\log t + x \right) \leq a_6e^{  - b_6 x},$$
 for all $(t,x)$ satisfying $\eqref{pair}$.
\end{lm1}

\begin{proof}
Note that $\sup_{u\leq t}|\Phi_1(u)|\leq \max_{k\leq m(t)+1}\eta_k,$ since $m(u)$ is the last renewal point of $T$ happening before $u,\,u\leq t.$ Consequently, for a fixed $c_6>0$ we have
$$\p\left(\sup_{u\leq t}|\Phi_1(u)|\geq c_6\log t + x\right)\leq \p\left(m(t) >\frac{2t}{\mu}\right) +\p\left(\max_{k\leq 2t/\mu+1}\eta_k > c_6\log t + x\right)=:\Phi_{1,1}+\Phi_{1,2}. $$
Using Lemma \ref{estimate_for_m} and assuming $a_3\geq 1,\,b_3\leq 1$  we get
$$\Phi_{1,1} \leq \p\left(m(t) - \frac{t}{\mu} - \frac{1}{\lambda\sqrt{\gamma}} \widetilde{W}_t >\frac{t}{2\mu}\right) + \p\left(\frac{\widetilde{W}_t}{\lambda\sqrt{\gamma}}>\frac{t}{2\mu}\right)  \leq a_3\exp\left\{c_3\log t - b_3\frac{t}{2\mu}\right\} + e^{-t/(8\gamma)} \leq $$ $$\leq
a_3\left(\frac{4\mu c_3}{b_3e}\right)^{c_3}e^{-b_3 t/(4\mu) } + e^{-x \lambda r/24}\leq a_3\left(\frac{4\mu c_3}{b_3e}\right)^{c_3}e^{-b_3 x  r/12 } + e^{-x \lambda r/24}. $$
Also,
$$\Phi_{1,2} \leq \frac{2t+\mu}{\mu}\exp\{-c_6 s_1 \log t  -s_1 x \}\e\exp\{s_1 \eta_1\}\leq \frac{2+\mu}{\mu}\e\exp\{s_1 \eta_1\}e^{-s_1 x} $$
 provided that $c_6$ is so large that $c_6 s_1 > 1.$
\end{proof}

To handle $\Phi_2$   we need a series of  auxiliary lemmas. The first of them is a modification of the Lemma 1 in Cs\"{o}rg\H{o} and Steinebach \cite{CS}.

\begin{lm1}\label{CsorgoSteinebach}
Suppose that $X_1,X_2,\ldots$ are i.i.d. centered random variables such that $\e \exp\{s X_1\}<\infty$ for some $s>0,$ and let $Q_n=\sum_{j=1}^n X_j$ $(Q_0=0).$ Then for any $L>0$ there exist $\delta = \delta(L),a_7,b_7,c_7>0$ $($which depend on $L,\delta$ and the distribution of $X_1)$ such that for any $n\in\mathbb{N}$ and all $x\geq 0$ one has
$$ \p\left(\max_{j\leq n}\max_{k\leq L\log n +\delta x} (Q_{j+k} - Q_j)\geq c_7 \log n + x\right) \leq  a_7e^{-b_7 x}.$$
\end{lm1}

\begin{proof} Clearly we can consider only non-constant distribution of $X_1.$ The probability under consideration is not greater than
$$(n+1)\sum_{j=1}^{[L\log n+\delta x]}\p(Q_j \geq c_7 \log n + x)\leq (n+1)\exp\{-sx - c_7s\log n\}\sum_{j=1}^{[L\log n+\delta x]}\e\exp\left\{ sQ_j \right\}\leq  $$  $$\leq (n+1)\exp\{-sx - c_7s\log n\}\frac{\varphi(s)}{\varphi(s)-1}\varphi(s)^{L\log n+\delta x}  \leq $$ $$\leq (n+1)\frac{\varphi(s)}{\varphi(s)-1} \exp\Bigl\{-sx - c_7s\log n + L(\log \varphi(s))\log n + (\log \varphi(s))\delta x  \Bigr\} , $$
where $\varphi(s)=\e \exp\{s X_1\}> 1.$ It   remains to take $c_7$ and $\delta$ so that inequalities  $c_7 s > L\log \varphi(s) + 1$ and $\delta(\log \varphi(s))<s$ are ensured.  
\end{proof}



 \begin{lm1}\label{Phi2}
There are  $a_8,b_8,c_8>0$ such that   $$\p\left(\sup_{u\leq t}  |\Phi_2(u)| \geq c_8\log t + 
  x \right) \leq a_8e^{  - b_8 x},$$
 for all $(t,x)$ satisfying $\eqref{pair}$.  
\end{lm1}

\begin{proof} Clearly we can consider only one coordinate of $S,$ say the first one. Also, using the remark before \eqref{greeks}, we can consider only $t\geq e\max\{1,\mu/2\}.$
Write, for $L>0$ and $\delta>0,$
$$\widetilde{L} := \frac{L}{2(1+\max\{0,\log \mu/2\})}, $$
$$\p\left( \sup_{u\leq t}|S_1(T_{m(u)}) - S_1(T_{[y]})|\geq  c_8\log t + x \right)\leq \p\left(N^{-1}(t/\gamma)>2t/\mu\right) + $$ $$+\p\left(\sup_{u\leq t}|m(u) - y(u)|> \widetilde{L}\log t + \delta x-1\right)+
\p\left(\sup_{\substack{n\leq 2t/\mu\\ k\leq \widetilde{L}\log t +\delta x}}  |S_1(T_{n+k}) - S_1(T_n)|\geq  c_8\log t + x\right)\leq $$ $$\leq\p\left(N^{-1}(t/\gamma)>2t/\mu\right)+
\p\left(\sup_{u\leq t}|m(u) - y(u)|> (\widetilde{L}-1)\log t + \delta x\right)+$$
$$+\p\left(\sup_{n\leq 2t/\mu}\sup_{k\leq  {L}\log [2t/\mu] +\delta x}  |S_1(T_{n+k}) - S_1(T_n)|\geq  c_8\log t + x\right) =:\sum_{q=1}^3  J_q,$$
where we have used that $2\log [z]\geq \log z$ if $z\geq e.$
Suppose that $L$ has been chosen to ensure  $\widetilde{L}\geq c_3+1 $ where $c_3$ is like in \eqref{m-y}, and that the numbers $\delta=\delta(L),a_7,b_7,c_7$ are provided by Lemma \ref{CsorgoSteinebach}, 
applied when $X_1$ has the  distribution of $\xi_{11}$. Then by Lemma \ref{CsorgoSteinebach} we see that 
$J_3 \leq a_7 e^{-b_7x}  $
provided that $c_{8}$ is so large that $c_{8}\log t \geq c_7\log  (2t/\mu)$   for all $t$ considered. 
To prove the statement of the lemma it remains to notice that by Lemma \ref{Poisson} and the relation \eqref{m-y} one has
$$J_1+J_2 \leq 2 e^{  - \lambda (1-2/e)   r x/3}+a_3e^{  - b_3 \delta x}. $$
\end{proof}

Consider now the case when $x $ lies in a domain of large deviations.

\begin{lm1}\label{LD} The relation $\eqref{1}$ holds for all $x > 3t /(\mu r).$ 
\end{lm1}

\begin{proof} First, by standard properties of Brownian motion one has
$$\p\left(\sup_{u\leq t}|W_u| > \frac{x}{2}\right)\leq 4de^{-3x/(8\mu r)} $$
for any $x$ provided that the pair $(x,t)$ satisfies Lemma's condition.
 Next, write
$$\p\left(\sup_{u\leq t}|S(u) - \varkappa u| > \frac{x}{2}\right)\leq \p\left(\sup_{k\leq m(t)}|S(T_k) - \varkappa T_k| > \frac{x}{6}\right) + $$ $$+\p\left(\sup_{u\leq t}|S(T_{m(u)}) - S(u)| > \frac{x}{6}\right) + \p\left(\sup_{u\leq t}|\varkappa||T_{m(u)} - u| > \frac{x}{6}\right) \leq $$ $$ \leq
\p\left(\sup_{k\leq m(t)}|S(T_k) - \varkappa T_k| > \frac{x}{6}\right) + \p\left(\sup_{k\leq m(t)+1}\eta_{k} > \frac{x}{6}\right) + \p\left(\sup_{k\leq m(t)+1}\tau_{k} > \frac{x}{6|\varkappa|}\right)\leq 
$$
$$\leq 3\p(m(t) > rx) + \p\left(\sup_{k\leq rx}|S(T_k) - \varkappa T_k| > \frac{x}{6}\right) + \p\left(\sup_{k\leq rx+1}\eta_k > \frac{x}{6}\right) + \p\left(\sup_{k\leq rx+1}\tau_k > \frac{x}{6|\varkappa|}\right) =\sum_{q=1}^4 R_q,
$$
By \eqref{deheuvels} and Lemma \ref{estimate_for_m}, assuming without loss of generality that $a_3\geq 1,b_3\leq 1,$ one has
$$\frac{R_1}{3} \leq \p\left(m(t) - \frac{t}{\mu} - \frac{1}{\lambda\sqrt{\gamma}}\widetilde{W}_t>\frac{rx}{2}\right)
+  \p\left( \frac{t}{\mu}  +\frac{\widetilde{W}_t}{\lambda \sqrt{\gamma}}>\frac{rx}{2}\right) \leq  a_3\exp\left\{c_3\log \frac{x\mu r}{3} - b_3\frac{rx}{2}\right\} + $$ $$+ \p\left( \frac{\widetilde{W}_t}{\lambda \sqrt{\gamma}}>r\frac{x}{6}\right)\leq  a_3\left(\frac{4\mu c_3}{3be}\right)^{c_3}e^{-b_3 rx/4} +  e^{-r^2x^2\lambda^2\gamma/(72  t)}\leq a_3\left(\frac{4\mu c_3}{3be}\right)^{c_3}e^{-b_3 rx/4} +  e^{-\lambda rx/24 } $$
for any $t\geq 1,x\geq 0$ if $x>3t/(\mu r).$

To estimate $R_2,$ by Doob's inequality for non-negative submartingales we get
$$R_2 \leq  \sum_{j=1}^d \p\left(\sup_{k\leq rx}\exp\{s_1(S_j(T_k) - \varkappa_j T_k)\} > e^{s_1 {x}/{6}}\right)
+ \sum_{j=1}^d\p\left(\sup_{k\leq rx}\exp\{s_1(\varkappa_j T_k - S_j(T_k) )\} > e^{s_1{x}/{6}}\right) 
 \leq $$ 
 \begin{equation}\label{R2}\leq  e^{-s_1   {x}/{6}}\sum_{j=1}^d\left(\e  \exp\{s_1(S_j(T_{[rx] }) - \varkappa_j T_{[rx] })\}  +   \e  \exp\{s_1( \varkappa_j T_{[rx] } - S_j(T_{[rx] }) )\} \right) \leq \end{equation} $$\leq 2d \exp\left\{-s_1  \frac{x}{6}      +   [rx] \log \e\exp\{s_1|\xi_1-\varkappa\tau_1|\} \right\}  \leq 2d e^{-qx}   $$
where $q:=s_1/6-r \log \e\exp\{s_1|\xi_1-\varkappa\tau_1|\}>0$ by \eqref{select_r}.
Now we estimate $R_{3}$ (and then $R_{4}$ is estimated analogously). From \eqref{select_r} it is seen that we can assume $r<s/12,$ in which case we have
$$R_{3}=\p\left(\sup_{k\leq rx+1}\eta_k > \frac{x}{6}\right)   \leq (rx+1) \p\left(\eta_1 > \frac{x}{6}\right) \leq (rx+1)e^{-s {x}/{6}}\e e^{s\eta_1}  \leq    \e e^{s\eta_1}e^{-sx/12}  . $$
\end{proof}

Theorem now follows from Lemmas \ref{Poisson}, \ref{MerlevedeRio}, \ref{Phi1}, \ref{Phi2}, \ref{LD}.

\end{proof}

We now turn to some corollaries of Theorem \ref{TheoremMain}. The first one concerns sums of a random amount of random summands.

\begin{cor1}
Consider a sequence of i.i.d. random vectors $\{(\xi_n = (\xi_{1n},\dots ,\xi_{dn}),\tau_n),n\geq 1\}$ with values in $\mathbb{R}^{d}\times\mathbb{R}_+$.  Let there exist $s>0$ such that $\e \exp\{s \tau_1\}<\infty$ and $\e\exp\{s |\xi_{k1}|\}<\infty$, $k=\overline{1,\, d}$. Denote $m(t) = \max\left(k\geq 0 : \sum\limits_{i=1}^k \tau_i\leq t\right)$ and
$S(t) = \sum\limits_{i=1}^{m(t)} \xi_i$ $($assuming a sum over an empty set is zero$)$. Then one can redefine the sequence $\{(\xi_n,\tau_n)\}$ on a probability space $(\Omega,\mathcal{F},\mathsf{P})$ supporting a standard Wiener process $\{W_t,t\geq 0\}$ such that for some positive constants $a,b,c,$ any $x>0$ and  $t\geq 1$ one has
\begin{equation}\label{RandomSumAsa}
\p\left(\sup_{u\leq t}\left|S(u) - \varkappa u - \sigma W_u\right| \geq c\log t + x \right) \leq a e^{-b x}, 
\end{equation}
where $\varkappa = {\e\xi_1}/{\e\tau_1}$ and $\sigma^2 = \di (\xi_1 - \varkappa\tau_1 )/\e\tau_1.$
\end{cor1}
\begin{proof}
It is sufficient to note that $S(u)$ is a cumulative process satisfying the conditions {\bf(A)} and {\bf(B)} with $T_0=0$ and $T_k = m^{-1}(k-1),\,k\geq 1$. Then, $\eta_k \equiv |\xi_k|$ and all the conditions of Theorem \ref{TheoremMain} are fulfilled.
\end{proof}

The second corollary  deals with birth and death processes. Recall that a continuous time homogeneous Markov chain $X=\{X_t,t\geq 0\}$ taking values in $\mathbb{Z}_+$
is an irreducible birth-death process if its transition matrices $\{P(t),t\geq 0\}$ satisfy 
$$\frac{dP(t)}{d t} = P(t)A = AP(t), $$
where the entries of infinitesimal matrix $A$ are
$$a_{i,i+1}=\lambda_i ,i\geq 0;\;a_{i,i-1}=\mu_i,i\geq 1;\; a_{0,0} = -\lambda_0,\, a_{i,i} = -\lambda_i-\mu_i,i\geq 1\mbox{ and }a_{i,j}=0\mbox{ if }  |i-j|>1. $$
Denote 
$$\pi_0=1,\;\pi_n = \frac{\lambda_0\ldots \lambda_{n-1}}{\mu_1\ldots \mu_n},\;n\geq 1. $$
Karlin and McGregor  \cite{KarlinMc} proved that conditions $\sum_n\pi_n<\infty,$ $\sum_n(\lambda_n\pi_n)^{-1}=\infty$ ensure that $X$ is ergodic, in which case $X$ satisfies the condition {\bf{(A)}}, with $T_0=0$ and $\{T_k,k\geq 1\}$ taken to be the moments of returning to zero.
For the condition {\bf{(B)}} to be true, we will require two following assumptions: 
\begin{equation}\label{vanDoorn}
\liminf_{n\to\infty} (\lambda_n+\mu_n - \sqrt{\lambda_{n-1}\mu_n} - \sqrt{\lambda_{n}\mu_{n+1}} )>0 
\end{equation}
and 
\begin{equation}\label{exponential}
\limsup_{n\to\infty}\left(\lambda_n\pi_n\right)^{1/n} <1.
\end{equation}
Note that \eqref{vanDoorn} and \eqref{exponential} together imply the Karlin-McGregor ergodicity conditions \cite{KarlinMc}, in which case there exists a stationary distribution $\tilde{\pi}$ on $\mathbb{Z}_+,$ defined as
$$\tilde{\pi}_n = \frac{\pi_n}{\sum_{k=0}^{\infty}\pi_k},\,n\geq 0.$$ Furthermore,
by Van Doorn \cite{VanDoorn} and Karlin-McGregor \cite{KarlinMc} the relation \eqref{vanDoorn}  ensures that  the Laplace transform of the return time to zero is an analytic function at the origin. Finally,
by \cite[\S I.12]{Chung} \eqref{exponential} means that the maximum of the birth and death process between two returns to zero has a finite exponential moment.
Thus we arrive at the following corollary.

\begin{cor1}
Consider a birth and death process $X$ satisfying \eqref{vanDoorn} and \eqref{exponential}, with initial distribution $\tilde{\pi},$ and let $f:\mathbb{R}\to \mathbb{R}$ be a   function growing not faster than a linear one. Then one can redefine $X$ on a probability space $(\Omega,\mathcal{F},\mathsf{P})$ supporting a standard Wiener process $\{W_t,t\geq 0\}$ such that for some positive constants $a,b,c,$ any $x>0$ and  $t\geq 1$ one has
\begin{equation}\label{BirthDeathProcess}
\p\left(\sup_{u\leq t}\left|\int_0^u f(X_s) ds - \varkappa_f u - \sigma_f W_u\right| \geq c\log t + x \right) \leq a e^{-b x},
\end{equation}
where $\varkappa_f = \sum_{n=0}^{\infty}f(n)\tilde{\pi}_n$ and $$\sigma_f^2 =
2\int_0^{\infty}\sum_{n,m=0}^{\infty}f(n)f(m)\tilde{\pi}_n(P_{nm}(s)-\tilde{\pi}_m)  ds.$$
\end{cor1}
\begin{proof}
Given the discussion above, it suffices to compute $\varkappa_f$ and $\sigma_f.$ However by central limit theorem for Markov chains, 
$\varkappa_f=\e f(X_0)$  and
$$\sigma_f^2 = \lim_{u\to\infty} \frac{1}{u}\int_{0}^u\int_0^u cov(f(X_t),f(X_s))ds \,dt =
2\int_{0}^{\infty} cov(f(X_0),f(X_s))ds. $$
\end{proof}

\end{document}